\documentclass[a4paper]{article}
\pdfoutput=1

\usepackage{IEK10} 
\usepackage{natbib}

\makeatother
\usepackage{amssymb}
\usepackage{xcolor}
\usepackage{mathtools, cuted}
\usepackage{tabularx}
\usepackage{eurosym}
\usepackage[utf8]{inputenc}

\usepackage{pstricks, pst-plot}	
\usepackage{graphicx} 
\usepackage{wrapfig} 
\usepackage[figuresright]{rotating}
\usepackage{lipsum}
\usepackage{placeins}
\usepackage{etoolbox}

\usepackage{amsmath}
\numberwithin{equation}{section} 
\usepackage{amssymb}
\usepackage{setspace}
\usepackage{enumerate}
\usepackage{enumitem}
\usepackage[english]{babel}
\usepackage{blindtext}
\usepackage{epsfig} 
\usepackage[utf8]{inputenc}			
\usepackage{ifpdf} 							
\usepackage{multicol}						
\usepackage{titletoc}						
\usepackage{setspace}						
\usepackage{ulem} 							
\usepackage{longtable}						
\usepackage{multirow,bigdelim}
\usepackage{booktabs} 						
\usepackage{parcolumns}						
\usepackage{tabularx} 						
\usepackage{threeparttable}					
\usepackage{graphicx}							
\usepackage{flafter}							
\usepackage{placeins}							
\usepackage{rotating}							
\usepackage{tikz}								
\usetikzlibrary{external}										
\usepackage{pgfplots}							
\usepackage{epstopdf}							
\pgfplotsset{compat=newest}						
\pgfplotsset{plot coordinates/math parser=false}
\usepackage{pdfpages}							
\graphicspath{{../02_Grafiken/}} 
\usepackage{float}

\renewcommand{\min}[1][]{
	\ifthenelse{\isempty{#1}}{\operatorname{min}}{\ensuremath{\underset{#1}{\text{min}\,}}}
}

\def\IEK10{
Institute of Energy and Climate Research -- Energy Systems Engineering (IEK-10),
  Forschungszentrum Jülich GmbH,
  Wilhelm-Johnen-Straße,
  52425 Jülich,
  Germany
}
\def\LTT{
  Institute of Technical Thermodynamics, 
  RWTH Aachen University,
  Schinkelstraße 8,
  52062 Aachen,
  Germany
}
\def\ETH{
    Energy \& Process Systems Engineering,
    Department of Mechanical and Process Engineering, 
    ETH Zürich,
  Tannenstrasse~3,
  8092 Zürich,
  Switzerland
}

\newcommand{\mytitle}{This is SpArta: Rigorous Optimization of Regionally Resolved Energy Systems by \underline{Sp}atial \underline{A}gg\underline{r}ega\underline{t}ion \underline{a}nd Decomposition}

\newcommand{\affil}{
  \begin{itemize}[leftmargin=3mm, itemsep=0mm]
    \item[$^a$]\LTT
    \item[$^b$]\ETH
    \item[$^c$]\IEK10
  \end{itemize}
}

\def\firstAuthor{Christiane Reinert}
\newcommand{\myauthor}{\firstAuthor$^{a}$, Benedikt Nilges$^{a}$, Nils Baumgärtner$^{a}$, and André Bardow$^{a,b,c,*}$}
\newcommand{\correspondingauthor}{A. Bardow, \ETH \newline E-mail: abardow@ethz.ch}
\author{\myauthor}

\usepackage[
  colorlinks,
  linkcolor=blue,
  citecolor=blue,
  urlcolor=blue,
  pdftitle={\mytitle},
  pdfauthor={\firstAuthor}
]{hyperref}
\usepackage[capitalise, nameinlink]{cleveref}
\crefname{table}{Tab.}{Tab.}

\newcommand{\setpgfexternalcounter}[1]{
  \makeatletter%
  \pgfkeysgetvalue{/tikz/external/figure name}\myexternalname
  \expandafter\gdef\csname c@tikzext@no@\myexternalname\endcsname{#1}%
  \makeatother
}

\begin{document}

  \thispagestyle{firststyle}

  \begin{center}
    \begin{large}
      \textbf{\mytitle}
    \end{large} \\
    \myauthor
  \end{center}

  \vspace{0.5cm}

  \begin{footnotesize}
    \affil
  \end{footnotesize}

  \vspace{0.5cm}

  \begin{abstract}
Energy systems with high shares of renewable energy are characterized by local variability and grid limitations. The synthesis of such energy systems, therefore, requires models with high spatial resolution. However, high spatial resolution increases the computational effort. Here, we present the SpArta method for rigorous optimization of regionally resolved energy systems by \underline{Sp}atial \underline{A}gg\underline{r}ega\underline{t}ion \underline{a}nd decomposition. SpArta significantly reduces computational effort while maintaining the full spatial resolution of sector-coupled energy systems. 

SpArta first reduces problem size by spatially aggregating the energy system using clustering. The aggregated problem is then relaxed and restricted to obtain a lower and an upper bound.
The spatial resolution is iteratively increased until the difference between upper and lower bound satisfies a predefined optimality gap. 

Finally, each cluster of the aggregated problem is redesigned at full resolution. For this purpose, SpArta decomposes the original synthesis problem into subproblems for each cluster. Combining the redesigned cluster solutions yields an optimal feasible solution of the full-scale problem within a predefined optimality gap. SpArta thus optimizes large-scale energy systems rigorously with significant reductions in computational effort. We apply SpArta to a case study of the sector-coupled German energy system, reducing the computational time by a factor of 7.5, compared to the optimization of the same problem at full spatial resolution. As SpArta shows a linear increase in computational time with problem size, SpArta enables computing larger problems allowing to resolve energy system designs with improved accuracy.
\end{abstract}

\vspace{0.5cm}


\vspace{0.75cm}

\begin{flushleft}
  \leavevmode{\parindent=15mm\indent}
  \textbf{Highlights:}
  \begin{itemize}[leftmargin=20mm]
    \item SpArta: rigorous design optimization for sector-coupled large-scale energy systems
    \item Near-optimal solution with known solution quality
    \item SpArta reduces computational effort while maintaining full spatial resolution
    \item SpArta reduces computational time by a factor 7.5 in German energy system model
    \item Linear increase in computational time with problem size
  \end{itemize}
\end{flushleft}
\vspace*{5mm}

  \newpage

\section{Introduction}\label{sec:Intro}
Planning energy systems with high shares of fluctuating renewable electricity generation requires optimization models with high spatial resolution to account for these fluctuations and local restrictions, such as grid limitations \citep{Keck.2022}. Local conditions such as wind speed and solar irradiation can play an important role and therefore need consideration by spatially-explicit modeling. Spatial resolution can be modeled by nodes, which reflect the location of each modeled site or region, and connections, which allow import and export  between nodes. For each node, the energy balance is closed, i.e., energy production and import match energy consumption and export. Currently, the spatial resolution in energy models ranges from one node per country represented \citep{Becker.2014, Kan.2021} over a more detailed regional resolution \citep{Bogdanov.2021}, to up to several thousand nodes for a single country \citep{Robinius.2017}. 

Choosing a spatial resolution sufficient to reflect the real system is both challenging and important: \cite{Heuberger.2020} show that a too coarse spatial resolution underestimates greenhouse gas emissions and energy storage needs. \cite{MarkoMimica.2022} further find that cost can vary significantly depending on the chosen spatial resolution. Furthermore, quantifying the spatial distribution of burdens and benefits of the energy transition can be achieved only by spatially explicit modeling \citep{Sasse.2020}. Additionally, models with high shares of renewable energy but without sufficient spatial resolution underestimate the importance of local constraints, such as grid limitations \citep{Rauner.2016}. \cite{Horsch.2017} therefore stress the importance of high spatial resolution for the joint optimization of energy conversion and grid extension. 

However, high spatial resolution is computationally challenging due to the increased problem size and, thus, computational effort in memory usage and computational time. 
In practice, the computationally feasible resolution of energy models is limited. Sector-coupled energy systems are increasingly employed to model decarbonization beyond the electricity sector. Sector-coupling introduces additional complexity and, therefore, massively raises computational time \citep{BehrangShirizadeh.2022}. Thus, practitioners commonly simplify models of sector-coupled energy systems. Common simplifications in energy modeling are temporal or spatial aggregation \citep{Lopion.2018}. Recently, \cite{Martinez.2021} reviewed and evaluated 34 mostly sector-coupled energy models that employ various clustering methods for spatial aggregation to reduce problem size. They stress that while spatially explicit modeling gained importance within the last years, high spatial resolution also leads to computational challenges. \cite{Frew.2016} discuss tradeoffs between spatial and temporal resolution as critical drivers of computational effort, showing that computational time and maximum memory usage increase quadratically with problem size for their linear model. \cite{Frew.2016} further list typical problem sizes of energy system models which are between 10 and 256 regions and range between 14 typical time slices up to an hourly resolution. The typical problem size for energy models also depends on other factors, such as the mathematical complexity and system scope \citep{Ridha.2020}. 

However, while aggregation is an effective tool to reduce problem size and computational effort, the solution resulting from the aggregated problem is usually infeasible for the full-scale problem. \cite{Frysztacki.2022} systematically evaluate the feasibility of spatially aggregated models and underline that the results of aggregated models can be infeasible at higher spatial resolution. \cite{Reinert.2020} and \cite{Syranidou.2020} thus discuss the necessity to disaggregate the resulting aggregated energy system and propose heuristic methods for disaggregation. \cite{Frysztacki.2022} show that optimization-based disaggregation leads to low load shedding compared to other approaches. However, heuristic methods are not able to guarantee system feasibility and cannot quantify the accuracy of the solution. In heuristic methods, the solution quality can only be assessed by comparison with the solution of the full-scale problem. However, avoiding the solution of the full-scale problem motivated the aggregation in the first place. When the full-scale problem cannot be solved, e.g., due to computational limitations, the accuracy cannot be quantified. Hence, known solution quality and feasibility without load shedding can only be guaranteed when the aggregated solution overestimates the required infrastructure of the full-scale energy system.  

This challenge is well-known from the time domain: Here, time series aggregation is often employed to aggregate the temporal resolution \citep{Teichgraeber.2022, Hoffmann.2022}, sometimes while taking into account the location dependence of the time series \citep{Weimann.2022}. Recent methods also exploit the ability of high-performance computers to parallelize problems \citep{Rehfeldt.2022}. Rigorous methods have been proposed that solve the full-scale problem to a specified accuracy without needing the full-scale problem as a benchmark: During the temporal aggregation, the objective function is under- and overestimated to rigorously solve optimization problems and obtain feasible solutions with a specified optimality gap \citep{Bahl.2018, Baumgartner.2019}. 

For the spatial dimension, \cite{Ma.2021} present such a rigorous approach for the aggregation of supply chains. However, their method is not applicable to energy systems with more complex electricity flow models and for products with limited transport opportunity, such as heat. Furthermore, the focus of their work is on identifying the near-optimal objective value. However, energy system design requires not only the objective value, but also a feasible full-scale design at the original spatial resolution. \cite{Moolya.2022} build up on our heuristic method for aggregation and disaggregation \citep{Reinert.2020} to rigorously optimize well placement in oilfields for multi-period optimization. However, their approach can not be directly transferred to energy system models which are strongly coupled and employ advanced transport models.
Therefore, a rigorous method to identify a feasible full-scale design for sector-coupled energy systems at large scale is still missing.

\subsection{Contribution}
In this work, we introduce SpArta, a method for the rigorous optimization of regionally resolved energy systems by spatial aggregation and decomposition. SpArta bridges the gap between solving sector-coupled energy systems with high spatial resolution and reduced computational effort. SpArta finds a rigorous, near-optimal solution at a specified accuracy which is feasible for the full-scale problem.
This work extends heuristic disaggregation methods by providing a guaranteed solution quality: To the best of the authors knowledge, SpArta is the first rigorous spatial aggregation method for optimizing sector-coupled energy systems with electricity transmission that ensures the solution’s feasibility in the full-scale problem. Optimality is guaranteed within a predefined optimality gap; feasibility is guaranteed by decomposing an aggregated solution to obtain a design and operation for the full-scale problem.  

We describe the SpArta method in Section \ref{sec:method}. In Section \ref{sec:casestudySpArta}, an example of a regionally resolved energy system with sector-coupling is given by optimizing the German energy system using SpArta. Section \ref{sec:conclusion} provides our conclusions.

\FloatBarrier
\newpage
\section[The SpArta method]{SpArta: A method for rigorous optimization by spatial aggregation and decomposition}\label{sec:method}
\begin{figure}[ht]
	\centering
	\includegraphics[width=1.0\textwidth]{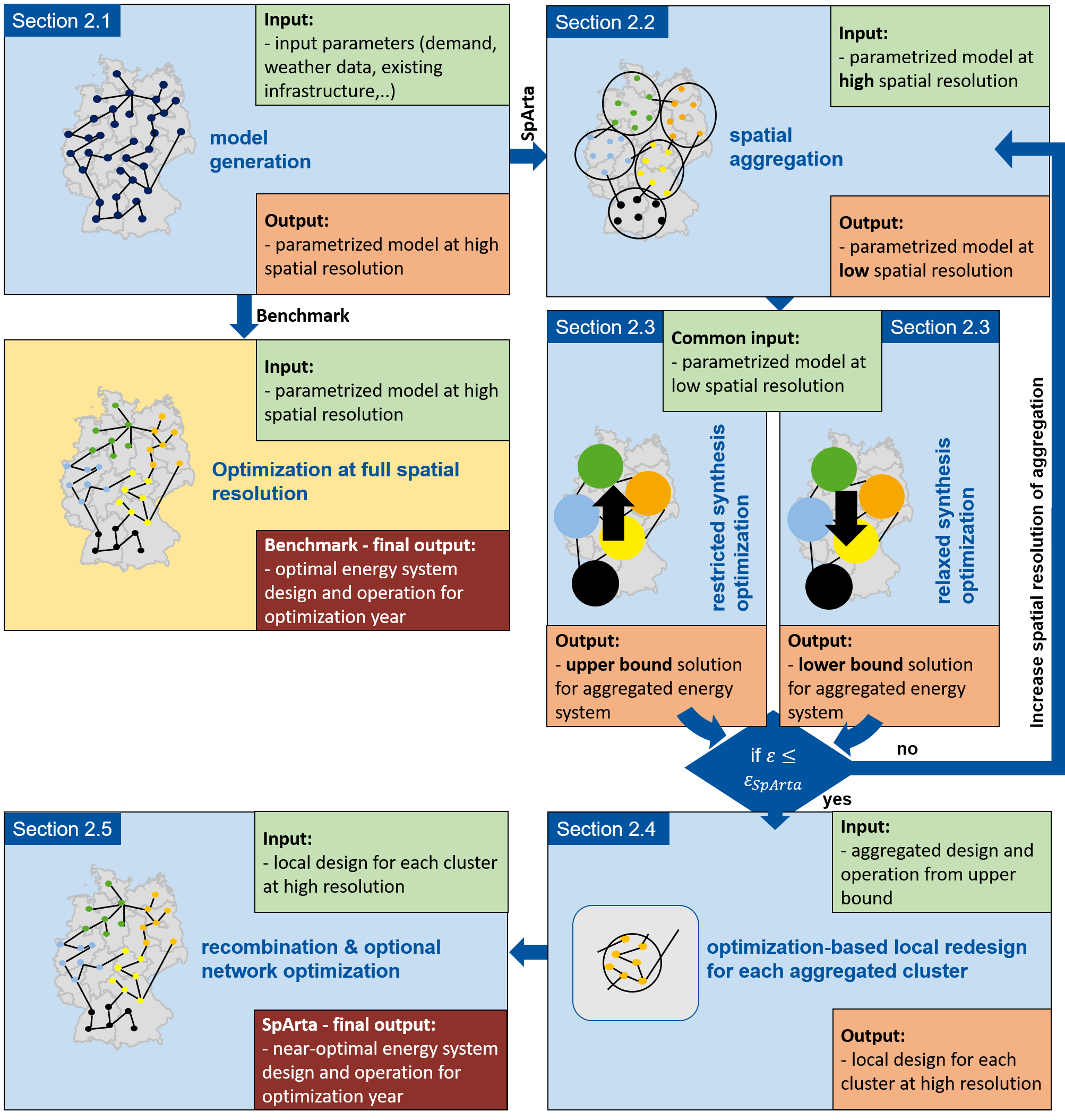}
	\caption{SpArta method for rigorous optimization of regionally resolved energy systems by spatial aggregation and decomposition: First, the regionally resolved input data is spatially aggregated. Then, lower and upper bounds are computed and compared to satisfy an optimality gap $\varepsilon_{SpArta}$. The spatial resolution of the clusters is increased until the optimality gap is reached. Then, the aggregated design is decomposed, and last, the resulting regionally resolved design is optimized in a network optimization. We optimize the system at full spatial resolution as a benchmark (yellow).}	
	\label{fig:SpArtaMethod}
\end{figure}

SpArta is a method to reduce computational time to find an optimal solution of highly resolved and sector-coupled energy systems within a specified optimality gap. In this Section, we first introduce a generic model formulation for energy system optimization (Section \ref{sub:GenericproblemDefinition}). Using this formulation, we discuss the steps of the SpArta method (Figure \ref{fig:SpArtaMethod}). 
In the first step, the SpArta method aggregates the system to clusters to simplify the synthesis problem by reducing spatial resolution (Section \ref{sub:Iterativespatialaggregation}). In the second step, we then both relax and restrict the spatial constraints in an aggregated optimization. The relaxation simplifies the problem to obtain a lower bound. In contrast, the restriction overestimates the problem to get an upper bound (Section \ref{sub:Lowerboundupperbound}). Iteratively, we solve both the relaxed and restricted synthesis problems at increasing spatial resolution until the difference between lower and upper bound solution values satisfies the specified optimality gap $\varepsilon_{SpArta}$. The upper bound overestimates the objective function value of a feasible solution for the full-scale problem. However, due to the aggregation, the upper bound solution does not distinguish where energy infrastructure is placed. If we want to obtain a feasible solution providing the detailed spatial design, we need a further step to resolve the infrastructure: Once the solution quality $\varepsilon$ satisfies the optimality gap $\varepsilon_{SpArta}$, we thus perform step three, where we consider each individual cluster but now at its full spatial resolution as subproblem. We optimize the subproblems as individual problems, using constraints from the results of the aggregated optimization, such as overall capacity in the cluster (Section \ref{sub:Decompositionandlocalredesign}). Last, we recombine the resulting design and operation of all subproblems to obtain an optimal solution for the full-scale problem. Advanced grid modeling approaches require a network optimization as additional step: A design optimization of the electricity grid at fixed energy converter infrastructure ensures feasibility at full scale. 
\FloatBarrier
\subsection{Generic problem definition} \label{sub:GenericproblemDefinition}
SpArta reduces the computational effort of the design optimization of multi-sectoral energy systems with high spatial resolution. In the following, we first introduce a generic model formulation. In general, an optimization problem can be formulated as minimization problem with equality constraints and inequality constraints \citep{Nocedal.1999}.

Our model formulation is based on \cite{Baumgartner.2021} with modifications to represent a more generalized system (as available as open-source framework in \cite{Reinert.2022}). The problem statement is as follows: 
\setlist{nosep}
\begin{itemize}[noitemsep]
	\item The system has an exogenous demand $d_{b,n,t}$ for each product $b \in B$ (e.g., electricity, heat, and transport) that is spatially resolved in nodes $n \in N$ and temporally resolved in time steps $t \in T$ of length $\Delta t_{t}$ in each year. Additionally, endogenous product demand may occur by the operation of components $c \in C$.
	\item A set of components $C^{prod}\subset C$ (e.g., wind turbines, heaters, and gas turbines) can be used for production $p$ of products $b\in B$ at each node and time step.
	\item Transport components $C^{grid} \subset C$ can move products between nodes via edges $l \in L$, of length $\Delta l_l$ (e.g., electrical grids).
	\item Greenhouse gas (GHG) emissions of the system are restricted.
\end{itemize}
The objective is to find the optimal energy system that minimizes the total annualized cost while obeying a GHG limit. The total annualized cost comprise the annualized investment cost (CAPEX) for all components built in any investment year $y \in Y$, including the current investment year $\textup{y}^*$ and previous investment years, and the operating cost (OPEX) in the current investment year. 

\begin{figure}[ht]
	\centering
	\includegraphics[width=\textwidth]{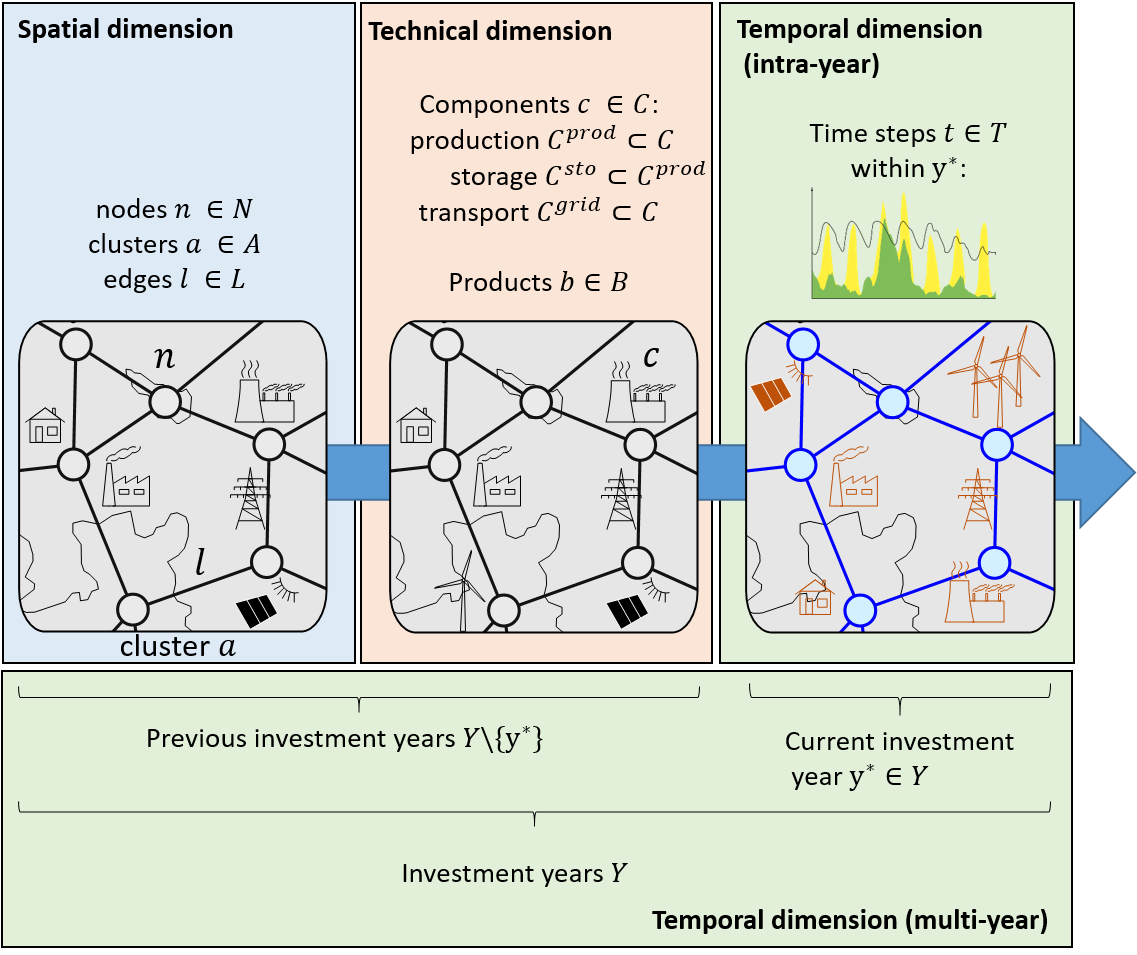}
	\caption{Spatial, temporal and technical dimensions in SpArta.}	
	\label{fig:Setfig}
\end{figure}

Figure \ref{fig:Setfig} and Table \ref{tab:Nomenclature} give a brief overview of all sets, superscripts, parameters, and variables used to describe the generic energy model and the SpArta method. Please note that for brevity, the production components $C^{prod}\subset C$ also include storage components, for which additional equations apply that are not discussed in detail in this work.

\subsubsection{Objective function and decision variables}\label{subsub:objective function}
The optimization problem minimizes the objective function while satisfying the demand $d$ for multiple products $b \in B$ in the current investment year $\textup{y}^*\in Y$ under environmental and technical constraints (Eq. \ref{eq:objective}). Here, we use the total annualized cost (TAC) of product supply as objective function, however, SpArta is generally not limited to TAC as objective function. Decision variables refer to infrastructure and operation. Infrastructure-related decision variables are the nominal production capacity expansion $p_{c,n,\textup{y}^*}^{nom}$ and the nominal transport capacity expansion $p_{c,l,\textup{y}^*}^{flow,nom}$ for the current investment year $\textup{y}^*\in Y$. Further, we have operational decisions $p_{c,n,t}$ for all components $c\in C$ and import decisions $p_{b,t}^{imp}$. All variables refer to a component-specific reference unit (e.g., $1$ kW\textsubscript{el}); respectively, the component-specific costs are also given per specific reference unit (definition ranges given in Table \ref{tab:Nomenclature}). 
\begin{equation}
\label{eq:objective}
\begin{split}
	\min_{p_{c, n, \textup{y}^*}^{nom}, p_{c, l, \textup{y}^*}^{flow, nom}, p_{c, n, t}, p_{b, t}^{imp}} \textup{   }
	&\underbrace{\sum_{c \in C^{prod}} \sum_{n \in N}\sum_{y\in Y} \frac{k_{c, y}^{inv, eco}}{\beta_c} p_{c, n, y}^{nom}}_{\text{CAPEX\textsuperscript{prod}}}\\
	+ &\underbrace{\sum_{c \in C^{grid}} \sum_{l \in L}\sum_{y\in Y} \frac{k_{c, y}^{inv, eco}}{\beta_c} \Delta l_l p_{c, l, y}^{flow, nom}}_{\text{CAPEX\textsuperscript{grid}}}\\
	+ &\underbrace{\sum_{c \in C^{prod}} \sum_{n \in N}\sum_{t\in T} k_{c}^{op, eco} p_{c, n, t} \Delta t_t+ \sum_{b\in B} \sum_{t\in T} k_{b,t}^{op, eco, imp} p_{b, t}^{imp} \Delta t_t} _{\text{OPEX}}
\end{split}
\end{equation}
The initially existing capacities of production (including storage) and transport components, which were built in previous investment years, are exogenous parameters of the investment optimization.
If a transition path is calculated, additional capacity can be added in each investment year. The existing capacity is then expanded by the capacity determined in the previous investment year and reduced by by decommissioning of components arriving at their end of lifetime.
The investment costs $k_{c,y}^{inv,eco}$ depend on the investment year $y$. The investment costs of all production components CAPEX\textsuperscript{prod} are defined as the overall capacity $p_{c,n,y}^{nom}$ multiplied by the annualized specific investment costs $\frac{k_{c,y}^{inv,eco}}{\beta_c}$ , where $\beta_c$  is the net present value factor \citep{Broverman.2017}, calculated with the interest rate $i$ for the time horizon $h_c^{min}$.
\begin{equation}
	\beta_c = \frac{(1+i)^{h_c^{min}}-1}{(1+i)^{h_c^{min}} i}
\end{equation}
As time horizon, we here use the minimum of the user-defined maximum discounting period $h_c$ and the actual component lifetime $h_c^{max}$.
\begin{equation}
	h_c^{min} = \min(h_c^{max}, h_c)
\end{equation}

Similar to the investment costs of all production components, we calculate the investment costs of all transport capacity (CAPEX\textsuperscript{grid}), multiplying the annualized specific costs per distance  $\frac{k_{c,y}^{inv,eco}}{\beta_c}$ with the distance $\Delta l_l$ and the overall transport capacity $p_{c,l,y}^{flow, nom}$. 

Operating costs for each time step $t$ depend on the length of the time step $\Delta t_t$, the specific operational costs $k_{c}^{op,eco}$ and used capacity $p_{c,n,t}$. Operating costs (OPEX) include maintenance and product imports. The operating costs for imports depend on specific import costs for each product $k_{b,t}^{op,eco,imp}$ and the overall product import $p_{b,t}^{imp}$.
\subsubsection{Product ratio matrix}
For each component, we define a fixed product ratio matrix $\theta_{b,c}$  that quantifies how much product is consumed, produced, charged/discharged, or transported if a reference unit of the component is used. As an example, for the production component gas turbine, the matrix entry $\theta_{b,c}$ denotes the ratio of natural gas consumed to produce a certain amount of electricity. As we discuss an LP problem, the product ratio matrix here denotes constant factors. For transport components $c \in C^{grid}$ , $\theta_{b,c}$ denotes which product can be transported via an edge when a reference unit of the component is installed.
\subsubsection{Supply-demand constraints}\label{subsub:supply-demand constraints}
The purpose of a multi-energy system is to satisfy the product demand at all time steps and nodes. All demands have to be matched by imports or production, including storage (Eq. \ref{demand_per_node}): The sum of the overall production $\theta_{b,c} p_{c,n,t}$ and product import $p_{b,t}^{imp}$ reduced by the internal losses $p_{b,t}^{flow,loss}$ must satisfy the exogenous nodal demand $d_{b,n,t}$. The transport losses depend on the transport efficiency $\eta_{c,l}^{grid}$.
\begin{equation}
\begin{split}
	& \sum_{c\in C^{prod}}\sum_{n \in N} \theta_{b,c} p_{c,n,t} + p_{b,t}^{imp}
	- \underbrace{ \sum_{l\in L} \sum_{c \in C^{grid}} \theta_{b,c} p_{c,l,t}^{flow}\left(1-\eta_{c,l}^{grid}\right) \Delta l_l}_{p_{b,t}^{flow,loss}} 
	\\
	 = &\sum_{n\in N} d_{b,n,t}\ ,\forall b \in B, \forall t \in T
\end{split}
\label{demand_per_node}
\end{equation}
The overall nodal production $\theta_{b,c} p_{c,n,t}$ may be positive (production) or negative (consumption). Using additional restrictions, we can limit the total imports $p_{b,t}^{imp}$ for some products $b$, for example, if a demand must be met without any imports $(p_{b,t}^{imp}=0)$. 

Production is limited by existing capacity from previous investment years and capacity expansion. Locally, for each component and node, the total production $p_{c,n,t}$  may not exceed the nominal capacity in any time step. In addition, the overall production may be limited by the temporal availability $\alpha_{c,n,t}$ (Eq. \ref{eq:availability}). The availability is especially crucial for fluctuating renewable energy converters.
\begin{equation}
\label{eq:availability}
	\alpha_{c,n,t} \sum_{y\in Y} p_{c,n,y}^{nom} \ge  p_{c,n,t} \ , \forall c \in C^{prod}, \forall n\in N, \forall t\in T
\end{equation}
In addition to the system-wide supply-demand constraints, the supply must match the demand at each node: the total production  $\theta_{b,c} p_{c,n,t}$ must be equal or exceed the sum of the exogenous nodal demand $d_{b,n,t}$ and the nodal exports $p_{b,n,t}^{exp}$ (Eq. \ref{eq:nodprod}). Nodal exports may be positive (export) or negative (import). In fluctuating energy systems, overproduction can occur when product supply exceeds the demand, e.g., for electricity. We allow overproduction at each node by an inequality constraint, to account for grid stability measures such as curtailment and, therefore, also to facilitate feasibility in the model.
\begin{equation}
\label{eq:nodprod}
	\sum_{c\in C^{prod}} \theta_{b,c} p_{c,n,t} - p_{b,n,t}^{exp} \ge d_{b,n,t} \ ,\forall b\in B, \forall n\in N, \forall t\in T
\end{equation}
In addition to the supply-demand constraints, we define an optional balance to meet a required minimum secured capacity $p_b^{min}$ for a product (Eq. \ref{eq:seccap}). 
\begin{equation}
\label{eq:seccap}
	\sum_{c \in C^{prod}} \sum_{n\in N} \sum_{y\in Y} \gamma_c \theta_{b,c} p_{c,n,y}^{nom} \ge p_b^{min} \ ,\forall b\in B
\end{equation}
Here, a capacity factor $\gamma_c$  denominates the share of a component’s capacity which is fully available at all times. In contrast to the time-dependent nodal availability $\alpha_{c,n,t}$, the capacity factor is an additional restriction to impose the overall long-term minimum availability of a component. The additional constraint can, therefore, e.g., reduce or avoid blackouts in the electricity sector for more detailed temporal resolutions or other weather years.

The overall capacity can be limited for some components, e.g., if capacity expansion is restricted for technological or political reasons. As an example, the overall wind converter capacity is usually limited. Hence, a capacity expansion is only possible if the sum over capacity expansion and existing capacity from previous optimization periods does not exceed a given limit $p_c^{max}$. The limit can be specified both at a node (Eq. \ref{eq:minnodcaplimit}) and for the full system (Eq. \ref{eq:mincaplimit}):
\begin{equation}
\label{eq:minnodcaplimit}
	\sum_{y\in Y} p_{c,n,y}^{nom} \le p_{c,n}^{max} \ , \forall c\in C, \forall n\in N
\end{equation}
\begin{equation}
\label{eq:mincaplimit}
	\sum_{y \in Y} \sum_{n \in N} p_{c,n,y}^{nom} \le p_c^{max} \ ,\forall c\in C
\end{equation}
\subsubsection{Transport constraints}
Transport components transport a product between nodes and thus spatially shift supply and demand. Two nodes are interconnected by an edge that enable the transport of products. Transport is limited by a maximum connection capacity for each time step (Eq. \ref{eq:capmaxtransport}).
\begin{equation}
\label{eq:capmaxtransport}
	- \sum_{y \in Y} p_{c,l,y}^{flow,nom} \leq p_{c,l,t}^{flow} \leq \sum_{y \in Y} p_{c,l,y}^{flow,nom} \ ,\forall c\in C^{grid}, \forall l\in L, \forall t\in T
\end{equation}
The product flow $p_{c,l,t}^{flow}$ is the flow between a node $n \in N$ and any other node $m \in N \backslash \{n\}$ connected by the edge $l$ . Generally, products can be transmitted in both directions. The connectivity matrix $\sigma_{l,n,m}$ contains all predefined connections of nodes and edges, i.e., two nodes are assigned to each connection. The overall export $p_{b,n,t}^{exp}$ for each node and product equals the product flow over all edges which connect the node with other nodes (Eq. \ref{eq:nodesedges}).
\begin{equation}
\label{eq:nodesedges}
	p_{b,n,t}^{exp} = \sum_{c\in C^{grid}} \sum_{l \in L}\theta_{b,c}  \sigma_{l,n,m} p_{c,l,t}^{flow} \ ,\forall b\in B, \forall n\in N, \forall m \in N\backslash \{n\}, \forall t\in T
\end{equation}
When a product is exported from $n \  (p_{b,n ,t}^{exp}>0)$, it is imported in $m \  (p_{b,m,t}^{exp}<0)$. $p_{b,n ,t}^{exp}$ denotes product exchange between nodes, whereas to the import variable $p_{b,t}^{imp}$ denotes imports over the system boundary. Optionally, an efficiency $\eta_{c,l}^{grid}$ can be defined for each component to account for transport losses.

If a product is transportable, the transport balances follow a simple transshipment approach for each product $b\in B$, except for electricity. For electricity, we use the DC-load-flow approach, as it is relatively common in large-scale energy models \citep{Overbye.2004}. The DC-load-flow approach is a simplification of electricity grids with alternating current and neglects losses. However, the DC-load-flow-approach shows high accuracy at much reduced computational time.

In the DC-load-flow approach, the nominal flow capacity $p_{c,l,y}^{flow,nom}$ is a property of the transmission line type, i.e., its voltage level, impedance, resistance, and the number of circuits. The electrical properties are used to calculate the susceptance $s_{c,l}$. $\Delta \vartheta_{c,l,t}$ describes the voltage difference between two nodes and $p_{c,l,t}^{flow}$ the electric product flow (Eq. \ref{eq:dclf}). 
\begin{equation}
\label{eq:dclf}
p_{c,l,t}^{flow} = s_{c,l} \Delta \vartheta_{c,l,t} \forall c \in C^{grid}
\end{equation}
The linear approximation by the DC-load-flow approach is described in depth by \cite{vandenBergh.2014}. Compared to more simple transport approaches such as transshipment, the DC-load-flow approach introduces additional constraints regarding the electricity flow. Generally, the SpArta method can be used with all linear transport models for electricity. 

Transport may not be applicable for all products, e.g., heat can usually not be transported. We therefore distinguish between transportable products, for which the transport components exist or can be built, and non-transportable products, which cannot be transported between nodes.
\subsubsection{Environmental constraints}
To account for environmental considerations in the design of energy systems, we here limit the annual operational greenhouse gas (GHG) emissions $e^{max}$ (Eq. \ref{eq:opimpactlimit}). 
\begin{equation}
\label{eq:opimpactlimit}
	e^{max} \ge \sum_{c \in C}\sum_{n \in N}\sum_{t\in T} k_{c}^{op,env} p_{c,n,t} \Delta t_t
\end{equation}
The operational impacts depend on the specific operational emissions $k_{c}^{op,env}$ for component operation. When accounting for operational emissions, we include the upstream emissions of products consumed during the conversion process, such as mining of fuels. However, the SpArta method can be generalized and is not limited to the specific GHG contraints.

In the following, we use the generic problem formulation to discuss the SpArta method. The SpArta method can be generally applied to LP (linear programming) and MILP (mixed-integer linear programming) problems, however, problem-specific modifications might be necessary. 
 
\subsection{Iterative spatial aggregation}\label{sub:Iterativespatialaggregation}
In SpArta, we under- and overestimate the full-scale problem at a reduced spatial resolution to reduce computational effort. We increase the number of spatial clusters in every iteration step $i$ until the solution quality $\epsilon_i$ satisfies the specified optimality gap $\epsilon_{SpArta}$ (Eq. \ref{eq:epssparta}). For each spatial resolution, we calculate lower and upper bound solutions for the total annualized cost (TAC) (see Section \ref{sub:GenericproblemDefinition}). The lower bound solution is a potentially infeasible relaxation, and the upper bound solution is a restriction that overestimates the objective value of a feasible solution. By using the total annualized cost of the lower bound as a denominator, we ensure that we overestimate the actual optimality gap.
\begin{equation}
\label{eq:epssparta}
	\epsilon = \frac{TAC_i^{UB}-TAC_i^{LB}}{TAC_i^{LB}} \le \epsilon_{SpArta}
\end{equation}
During the spatial aggregation, we allocate each node to exactly one cluster. The initial spatial resolution is set to two clusters unless stated otherwise by the user. We spatially aggregate using common clustering algorithms, e.g., the k-means \citep{Macqueen67somemethods}, k-medoids \citep{Vinod.1969}, or hierarchical algorithm \citep{Ward.1963}.

\begin{figure}[ht]%
\centering
\includegraphics[width=\textwidth]{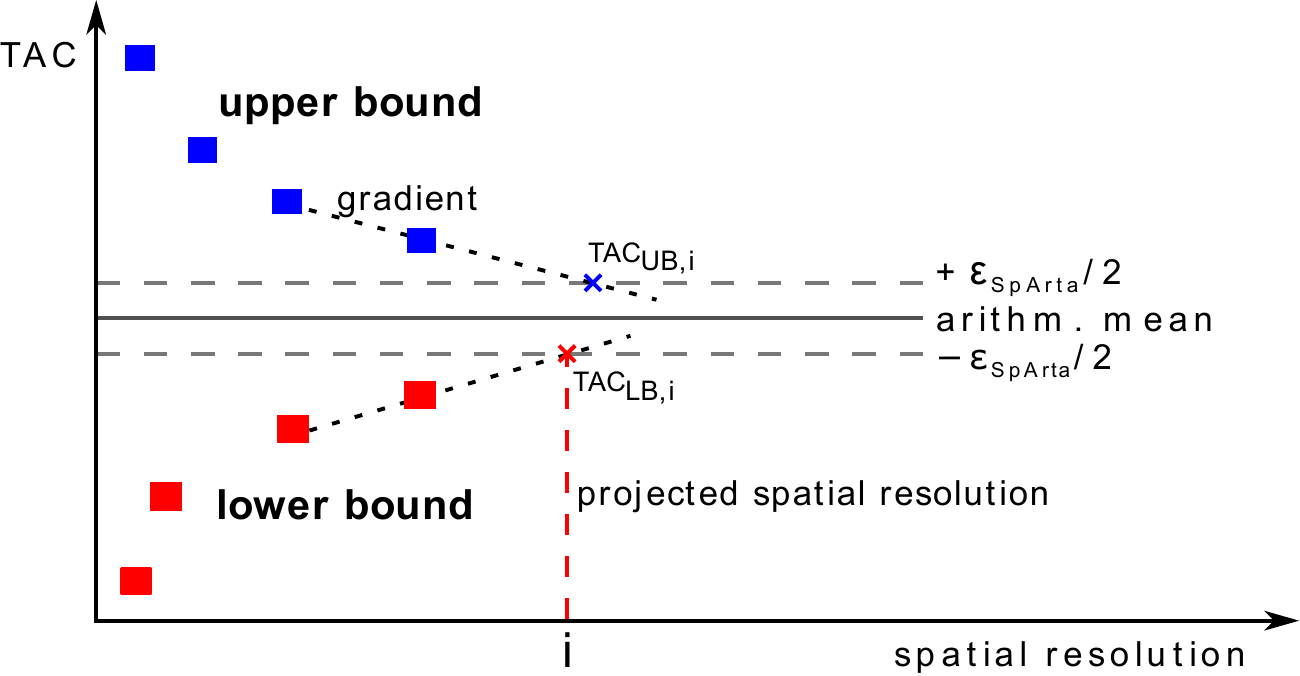}%
\caption[Fast-forward clustering algorithm]{Fast-forward clustering algorithm to determine the necessary spatial resolution till the predefined optimality gap is reached.}%
\label{fig:SpArta fast forward algorithm}%
\end{figure}

The spatial resolution is iteratively increased following a given rule, for example, a predefined increase in cluster elements. In this work, we further propose the fast-forward cluster algorithm to enhance the convergence towards the optimality gap. The fast-forward cluster algorithm increases the resolution as follows: We calculate gradients between two consecutive iteration steps for both the lower and upper bound.  We then extrapolate the gradients to estimate the expected required number of clusters at which the solution quality satisfies the predefined optimality gap (Figure \ref{fig:SpArta fast forward algorithm}): We calculate the intersection between the extrapolated gradient and the arithmetic mean of lower and upper bound plus/minus half the predefined optimality gap $\epsilon_{SpArta}$ for both the lower and upper bound. Estimating the necessary resolution of the spatial aggregation leads to faster convergence towards the optimality gap, compared to a continuous increase in cluster steps. 
The intersection occurs at a certain cluster number for both the lower and upper bound, resulting in two potential spatial resolutions. The smaller number of clusters is set as the spatial resolution of the next iteration. To further enhance convergence, the maximal and minimal increase of spatial clusters can be limited.
\FloatBarrier
\subsection{Lower \& upper bound} \label{sub:Lowerboundupperbound}
To rigorously under- and overestimate the optimization problem, we compute the lower and upper bounds of the full-scale problem from the problem at reduced spatial resolution. We relax constraints to obtain a lower-bound solution, thereby underestimating the objective value. This relaxation leads to possibly infeasible solutions that are not part of the original solution space, e.g., if an energy demand is underestimated, the resulting system design is usually not able to fulfill the demand of the original system. 

For the upper bound, we overestimate the restrictions of all limitations. The objective value of the aggregated problem exceeds the objective value of a feasible full-scale solution. Some results from the upper bound optimization are later used to find a feasible solution to the full-scale problem. In general, we still solve the problem from Section \ref{sub:GenericproblemDefinition} for the lower and upper bound in the aggregated system but with lower spatial resolution. Hence, unless stated otherwise, all balances for the whole system, such as the overall product balance, remain valid. Please also refer to the Nomenclature \ref{app:Nomenclature} for further detail on the used variables and parameters.

In the following, we discuss the restrictions and relaxations used to obtain lower and upper bounds. First, we discuss node-specific parameters (Section \ref{subsub:Node-specific parameters}). Then, we consider transport over cluster-internal connections (Section \ref{subsub:Internal connections}) and non-transportable products (Section \ref{subsub:Non-transportable goods}). 
\subsubsection{Node-specific parameters} \label{subsub:Node-specific parameters}
First, we restrict and relax the supply-demand constraints by aggregating the node-specific parameters, such as the nodal demand and existing capacities. For each cluster $a \in A$, we have a subset $N_a \subset N$, which comprises the nodes within the cluster $(n\in N_a)$. All edges within a cluster $a$ are defined as $l\in L_a^{int}$, whereas all edges that connect the cluster with other clusters are defined as $l\in L_a^{ext}$. $L_a^{int}\subset L$ and $L_a^{ext}\subset L$ are subsets of the overall edges $L$ in the full-scale problem. In the following, all aggregated variables and parameters for which the index $n$ is modified to the index $a$, are marked with a hat. As an example, $p_{b,c,n,y}^{nom}$ denotes the nominal capacity at node $n$, and $\hat{p}_{b,c,a,y}^{nom}$ denotes the nominal capacity at cluster $a$.\newline
\newline

For each cluster, we aggregate the existing infrastructure (built in previous investment years $y \in Y \backslash \{\textup{y}^*\}$) and exogenous product demand:
\begin{subequations}
\begin{equation}
\hat{p}_{b,c,a,y}^{nom} = \sum_{n\in N_a} p_{b,c,n,y}^{nom} \ ,\forall b\in B, \forall c\in C, \forall a \in A, \forall y \in Y \backslash \{\textup{y}^*\}
\label{eq:31a}
\end{equation}
\begin{equation}
\hat{d}_{b,a,t} = \sum_{n\in N_a} d_{b,n,t} \ ,\forall b\in B, \forall a\in A, \forall t\in T
\label{eq:31b}
\end{equation}
\end{subequations}
The overall nominal capacity $\hat{p}_{b,c,a,y}^{nom}$ of the previous investment years and the demand $\hat{d}_{b,a,t}$ are model-exogenous parameters. For the lower bound, we further overestimate the capacity availability $\widehat{av}_{c,a,t}$  for cluster $a\in A$ (Eq. \ref{eq:32}) by selecting the cluster-internal highest availability for each time step. 
\begin{equation}
 \widehat{av}_{c,a,t} = \max_{n\in N_a} \left( \alpha_{c,n,t} \right) \ ,\forall c\in C^{prod}, \forall a\in A, \forall t \in T
\label{eq:32}
\end{equation}
Respectively, we underestimate the availability for the upper bound (Eq. \ref{eq:33}) for all components $c\in C^{prod}$ and each cluster $a\in A$ by selecting the cluster-internal lowest availability for each time step.
\begin{equation}
\widehat{av}_{c,a,t} =\min_{n\in N_a} \left(\alpha_{c,n,t}\right) \ , \forall c\in C^{prod}, \forall a\in A, \forall t\in T
\label{eq:33}
\end{equation}

We perform two optimizations, one for the lower and one for the upper bound. In each optimization, we calculate the cost and environmental impacts according to the objective function (Eq. \ref{eq:objective}), replacing the nodes with clusters and the connections between the nodes with our estimations from Sections \ref{subsub:Internal connections} and \ref{subsub:Non-transportable goods}. 
\subsubsection{Transportable products} \label{subsub:Internal connections}
Within a cluster, we do not have any information about internal transport limitations. If nodes within a cluster are not connected, products cannot be transported between them - hence, production capacity at one node can not be used to meet the demand at another node. Sufficient product transport should be ensured to find tight bounds, as otherwise, the capacities within a cluster would be greatly overestimated. For this purpose, we extend the iterative spatial aggregation (Section \ref{sub:Iterativespatialaggregation}): The objective value of the upper bound can be reduced if we can assure that transport between all nodes within a cluster is possible for transportable products.  For this reason, we first check for each cluster whether all nodes are interconnected within the cluster. If the nodes within cluster $a\in A$ are not interconnected, we separate cluster $a$ into individual subclusters. Although the approach results in an overall higher number of clusters, computational experience shows benefits due to a better approximation of the problem. 

In the following, we under- and overestimate transport losses and capacity expansion of each cluster's internal connections to obtain a lower and upper bound. We use the under- and overestimated losses and capacity expansion in the supply-demand constraints (Section \ref{subsub:supply-demand constraints}) and the cost of transport infrastructure in the objective function (Section \ref{subsub:objective function}).
For the lower bound, we underestimate potential transport losses within a cluster by assuming ideal transport (Eq. \ref{eq:34}). 
\begin{equation}
\eta_{c,l} = 1 \ ,\forall c\in C^{grid}, l\in L_a^{int}
\label{eq:34}
\end{equation}
When ideal transport is assumed, the location of production and consumption within a cluster is not relevant for calculating the lower bound. Thereby, a product demand located anywhere in the cluster can be matched by production at any other cluster node without losses due to transport efficiency.

For the upper bound, we respectively overestimate potential transport losses. To obtain a worst-case scenario, we assume that the whole cluster's demand occurs at one single node. This node is assumed to be connected by the smallest existing internal transport capacity and has the lowest possible existing production capacity of all nodes in the cluster. The transport infrastructure within the cluster must be extended to fulfill the demand, endogenous consumption, and exports of any transportable product at this node.  To calculate the necessary capacity expansion, we first calculate the maximum product flow within a cluster $a$ as the sum of all demand and exports to all other clusters (Eq. \ref{eq:35}).

\begin{equation}
\hat{p}_{b,a,t}^{flow,max} = demand_{b,a,t} + export_{b,a,t} \ ,\forall b\in B, \forall a\in A,\forall t\in T
\label{eq:35}
\end{equation}
The demand $demand_{b,a,t}$ consists of both model-exogenous product demand $\hat{d}_{b,a,t}$ and model-endogenous product consumption. Product consumption is defined as negative production, hence denoted with a negative prefix:
\begin{equation}
demand_{b,a,t} = \hat{d}_{b,a,t} - \sum_{c\in C} \min \left(0, \theta_{b,c} \hat{p}_{c,a,t}\right).
\label{eq:36}
\end{equation}
Export from a cluster leads to an additional demand within the cluster. $export_{b,a,t}$ denominates the overall exports from cluster $a\in A$ to any other cluster $z\in A\backslash \{a\}$. As imports do not increase demand, we only consider positive exports:
\begin{equation}
export_{b,a,t} = \sum_{c\in C^{grid}}\sum_{z\in A\backslash\{a\}} \sum_{l\in L_a^{ext}} \max \left(\theta_{b,c}  \sigma_{l,a,z} p_{c,l,t}^{flow},0\right).
\label{eq:37}
\end{equation}
Using the maximum product flow, we overestimate cluster-internal transport losses (Eq. \ref{eq:38}) for the upper bound.
\begin{equation}
\hat{p}_{b,a,t}^{flow,loss,UB} = \left(1 - \prod_{l\in L_a^{int}} \min_{c\in C^{grid}}(\eta_{c,l})  \Delta l_l \right) \hat{p}_{b,a,t}^{flow,max} ,\forall b\in B,\forall a\in A,\forall t\in T
\label{eq:38}
\end{equation}
We assume that the maximum product flow must be transported over all transport elements within the cluster at each time step. Thus, as a worst-case assumption, we assume that the maximum flow is affected by the product of all transport efficiencies of all edges in the cluster. 

Further, we dimension the transport capacity of all cluster-internal edges to transport the maximum product flow to satisfy the demand. With the maximum product flow over all time steps, we can find an upper bound for the necessary additional nominal transport capacity $\hat{p}_{c,l,\textup{y}^*}^{flow,nom}$ within each cluster (Eq. \ref{eq:39}). 
\begin{equation}
\begin{split}
\theta_{b,c} \hat{p}_{c,l,\textup{y}^*}^{flow,nom} =& \max \left(0,  \max_{t\in T} \left(\hat{p}_{b,a,t}^{flow,max}\right) - \underbrace{\sum_{c\in C^{grid}} \sum_{y\in Y\backslash \{\textup{y}^*\}} \theta_{b,c} \hat{p}_{c,l,y}^{flow,nom} }_{\text{existing transport capacity}}\right) \textup{,} \\
& \forall c\in C^{grid},\forall b\in B,\forall l\in L_a^{int}, \forall a\in A
\end{split}
\label{eq:39}
\end{equation}

If the existing transport infrastructure at each cluster-internal edge is not sufficient to transport the maximum product flow, the transport capacity must be extended. In some cases, transport infrastructure cannot be extended due to capacity expansion restrictions. In this case, the maximum product demand at a node must be partly fulfilled by a local capacity expansion of production components. 

Cluster-external edges still exist in the spatially aggregated problem and therefore do not generally need any over- or underestimation of their capacity. However, for a more advanced transport model, such as the DC-load-flow approach, we need to ensure that during decomposition, edges cannot inhibit the power flow on other edges due to their additional load flow constraints. For the upper bound, we, therefore, overestimate the transmission capacity of all cluster-external edges $l\in L_a^{ext}$ which are modeled according to a DC-load-flow approach by the overall capacity of all cluster-external edges that connect the same cluster. The overestimation of the edges leads to a feasible upper bound, as the transmission capacity and, therefore, the overall cost is higher than in the unconstrained capacity expansion optimization where edges can be expanded up to their optimal capacity. 
\subsubsection{Non-transportable products} \label{subsub:Non-transportable goods}
In an aggregated model, we lose information on the spatial limitations of products for which no transport is possible. Hence, the aggregated system underestimates the full-scale problem: On one hand, the operating cost may be underestimated, as a component with low operating costs can satisfy product demands at nodes where only components with higher operating costs are installed. On the other hand, capacity expansion in the full-scale model may be underestimated in the aggregated system, as existing components at one node can be used to satisfy both product demands and the secured capacity constraint (as defined in Eq. \ref{eq:seccap}) at any node in the cluster. 

In the following, we employ node-specific parameters of the original problem, given for single nodes $n\in N_a$ in a cluster, to obtain cluster-specific constraints for non-transportable products.

To determine the operating cost of the aggregated system for the upper bound, we require an additional constraint, which must ensure that components are not used across nodes to fulfill product demands at different nodes over time. For each node and time step, we determine which existing components could be used to satisfy the exogenous demand at the least cost in advance of the optimization. To identify which components can be used in the restricted optimization, we first sort the operating costs of all components in a merit order. In each time step, any component operated in the aggregated system must, therefore, only be operated to the extent to which the product can be used to match the demand at the node in which it is located more cost-efficient than any other component at this node.

To ensure that newly built capacities are only operated to the extent that they could be operated in the fully resolved system, we must ensure that newly built capacity at one node cannot be used across nodes over time. If new capacity is installed in the aggregated problem, such new capacity could also be employed at any node inside the cluster to fill the demand. To avoid that the need for new components is underestimated as new capacity is used at more than one node over time, we ensure that existing components can only be operated when their operation is cheaper than the operation of any component that could be built. Any demand that cannot be satisfied by cost-efficient existing components then needs to be satisfied by sufficient capacity expansion at each node, which we ensure by an additional constraint (as later discussed in Eq. \ref{eq:regconstr}).

The merit-order strategy ensures that all existing components, which help satisfy the demand, can compete with capacity expansion: We determine the cost of a reference component by identifying the component with the lowest operating cost. We then sort the operating costs of all existing components in a merit order. If a component's operation is more expensive than the operation of the reference component, we do not allow its operation during the optimization. For each node and time step, we allow the use of existing components dependent on the merit order until the demand is satisfied or no cost-efficient components exist to further contribute to the production.   

For the lower bound, neglecting spatial restrictions generally poses a valid relaxation, as the aggregated solution yields a lower objective value than the full-scale problem. However, for the required secured minimal capacity (as defined in Eq. \ref{eq:seccap}), we introduce an additional constraint for both the relaxed and restricted optimization to improve the estimation (Eq. \ref{eq:41}).  

For the subset of all non-transportable products $B^{nt}\subset B$, the nominal capacity in the cluster is extended to match the needed secure capacity $\delta_{b,n}$, defined as the difference between the required minimal nodal capacity $p_{b,n}^{min}$ and the existing capacity that can be used at any time:
\begin{equation}
\delta_{b,n} = p_{b,n}^{min} - \sum_{c\in C^{prod}} \sum_{y\in Y\backslash \{\textup{y}^*\}} \gamma_c \theta_{b,c} p_{c,n,y}^{nom} ,\forall b\in B^{nt}, \forall n\in N_a .
\label{eq:42}
\end{equation}
If the existing infrastructure at a node $n$ is not able to fulfill the secured minimal capacity, additional infrastructure must be installed. The additional constraint ensures that an existing component in one node cannot be used to satisfy the minimal secured capacity constraint at another node in the same cluster. If the existing capacity in the cluster already exceeds the minimal required capacity in each node, no capacity expansion is needed. 
\begin{equation}
\sum_{n\in N_a} \left(\max\left(\delta_{b,n},0\right)\right) \le \sum_{c\in C^{prod}} \gamma_c  \theta_{b,c} \hat{p}_{c,a,\textup{y}^*}^{nom} \ ,\forall b\in B^{nt}, \forall a\in A
\label{eq:41}
\end{equation}

For the upper bound, we need a further estimation for the required capacity expansion. Besides the restriction for the secured capacity (Eq. \ref{eq:41}), we define the overall demand needs $\lambda_{b,n}$, as the maximum product demand over time at every single node that cannot be satisfied by the cost-efficient existing infrastructure at the node itself (see Eq. \ref{eq:44}). $\omega_{c,n,t}$ is a parameter between 0 and 1 and determines the useable cost-efficient existing capacity at each node for each component and time step.

Additional to the exogenous demand $d_{b,n,t}$, we consider the maximum model-endogenous consumption of each non-transportable product $\max \left(- \theta_{b,c}\hat{p}_{c,a,t},0\right)$ within a cluster. As we do not know at which node in the cluster the endogenous demand will occur, we introduce the number of all nodes within a cluster, the cardinality $\vert N_a\vert$, and assume this demand occurs at every single node within the cluster to rigorously overestimate the problem. Thereby, the $\lambda_{b,n}$ describes a pessimistic estimate of the necessary additional capacity needed to be able to fulfill the product demand, ensuring a sufficient capacity expansion in the aggregated system compared to the full-scale system. However, we would like to point out that depending on the problem, this constraint could likely be chosen stricter to obtain an improved bound when high non-transportable, model-endogenous demands occur.

\begin{equation}
\begin{split}
&\lambda_{b,n}
= \max_{t\in T} \left(d_{b,n,t} + \max \left(-\theta_{b,c}\hat{p}_{c,a,t},0\right) \vert N_a\vert - \sum_{c\in C^{prod}} \sum_{y\in Y\backslash \{\textup{y}^*\}} \theta_{b,c} \omega_{c,n,t} p_{c,n,y}^{nom},0\right) \\& ,\forall b\in B^{nt}, \forall n\in N_a
\end{split}
\label{eq:44}
\end{equation}

To find a rigorous upper bound for the capacity expansion of non-transportable products, we calculate the necessary capacity expansion to satisfy both the minimal required capacity $\delta_{b,n}$ and the overall product demand $\lambda_{b,n}$ at each node within the cluster:
\begin{equation}
\begin{split}
&\sum_{n\in N_a} \max\left(\delta_{b,n},\lambda_{b,n}\right)\\
&\le \sum_{c\in C^{prod}} \min_{t \in T}\left(\alpha_{c,a,t},\gamma_c\right) \theta_{b,c} \hat{p}_{c,a,\textup{y}^*}^{nom} \ ,\forall b\in B^{nt}, \forall a\in A.
\end{split}
\label{eq:regconstr}
\end{equation}

\subsection{Decomposition and local redesign} \label{sub:Decompositionandlocalredesign}
When the optimality gap is reached, we can bound the objective function of the full-scale problem. However, typically we are not only interested in the value of the objective function, such as the total annualized cost of a system, but also in the infrastructure design, i.e., the location of newly installed components.

To obtain a regionally resolved component distribution and thus a feasible solution for the full-scale problem, we decompose the aggregated problem into the single cluster nodes after reaching the optimality gap. Each cluster then consists of the nodes within one cluster of the aggregated problem during the final iteration. In the cluster, we use the upper-bound results of the aggregated problem as input parameters to redesign the energy system at full resolution: The overall nominal capacity of components in the cluster and the cross-border flows are set as fixed constraints to the synthesis optimization of each isolated cluster. Thus, the optimization of every single cluster locally redistributes the expanded component capacity resulting from the aggregated system to the full spatial resolution. We further optimize transport capacity expansion within each cluster. As an environmental constraint for each cluster, we limit the greenhouse gas emissions according to the cluster-specific greenhouse gas emissions for the upper bound solution. 

Combining the clusters' solutions yields a feasible design at full scale. As each cluster is redesigned independently, we reduce the degrees of freedom and thus the computational effort compared to the fully interconnected problem. The solution quality is equal to or exceeds the solution quality defined in the optimality gap: Due to fewer overestimations, the redesigned solution improves compared to the upper bound, and hence the solution quality is higher. 

\subsection{Operational and network optimization}\label{subsub:networkopt}
As the last step, the regionally resolved problem can further be optimized by an operational or network optimization to improve the solution quality:
\begin{itemize}
	\item {Operational optimization: An operational optimization of the full-scale problem with fixed transport and converter capacity can exploit synergies between the clusters and prove the feasibility of the full system. }
	\item {Network optimization: The operational optimization can further be extended by a network optimization to identify the optimal transport capacity distribution in the full-scale problem. While the energy converter design remains fixed, we allow capacity expansion of the transport components to further improve transport between the clusters in the full-scale problem.}
\end{itemize}
Despite the full spatial resolution, the computational time of the operational or network optimization is significantly smaller than the computational time of the full-scale problem, as the amount of variables is significantly lower due to the fixed energy converter designs.

With a simple transport approach, such as transshipment, the SpArta method leads to feasible solutions without a network optimization. If the DC-load-flow approach is used for electricity transport, the guaranteed feasibility of SpArta requires an additional network optimization. In the full-scale problem, if the power flow over an edge connecting two clusters reaches its maximum, the phase difference $\Delta \vartheta_{c,l,t}$ of this line becomes fixed and thus reduces the degree of freedom. This limitation may lead to infeasibility: If more than one phase difference is fixed, the phase difference of other edges at the shared node may be impacted as well. For example, if an edge within a cluster is neighbored by two edges that are operated at their power limit, the phase differences of the neighboring edges will determine the phase difference of the edge within the cluster. 

Thus, as a last step of SpArta, to ensure feasibility, a network optimization, i.e., a design optimization of the full grid, must be performed, such that potentially congested lines may be further extended. During the network optimization, the capacity of all production components is already determined and only considered as a parameter, leading to relatively low computational times. Even though feasibility cannot be guaranteed without the subsequent network optimization at full scale, in practice, we observed that this step could often be avoided, as the system is mostly already feasible. Alternative to directly applying a network optimization, feasibility can be tested by an operational optimization at full scale. 
\FloatBarrier

\section[Case study and results]{Case study and results: enhancing computational time by the SpArta method}\label{sec:casestudySpArta}

In this section, we briefly describe our case study, an energy system optimization model for Germany at high spatial resolution. We optimize Germany’s future energy system for the year 2030, comprising the electricity, heating, and transportation sectors. We consider the energy demand of both industrial users and private households in the electricity and heating sector, and private transportation in the transportation sector. Table \ref{secmod:techs} provides an overview of the technologies used to supply electricity, heating and transportation. In a brownfield optimization, we optimize design and operation of the energy system with regard to previously existing infrastructure, as given in \cite{Baumgartner.2021}.  We use the same input parameters for our case study and the benchmark. We obtain the design and operation of energy converters and grid infrastructure in a least-cost optimization. The optimization can expand the previously existing infrastructure, but is limited by local capacity expansion potentials. We simultaneously optimize the design and operation for all components, nodes, and time steps during the optimization.

The model accounts for regional fluctuations of energy demand and supply at a spatial resolution of 416 nodes, connected by 660 interconnections. Transport of electricity between the nodes is possible within the limits of grid infrastructure. We further consider non-transportable product demands, such as industrial and domestic heat demands. Further, we limited the temporal availability of components and the overall greenhouse gas emissions. The grid is modeled using a DC-load-flow approach (based on \cite{Egerer.2016}). Thus, transmission losses are neglected. Extension of the grid is possible without limitations by switching to a higher voltage level or by building new power lines to extend the existing grid.

\begin{table}[!ht]
    \centering
		\caption{Technologies used in our case study to supply electricity, heating, and transport.}
    \begin{tabular}{|p{0.5\textwidth}|p{0.5\textwidth}|}
    \hline
 \textbf{electricity} & \textbf{heating} \\ \hline
        biomass plants & thermal insulation (demand reduction) \\ 
        combined-cycle plants & electric boilers \\ 
        coal plants & gas boilers \\
        gas plants & oil boilers \\ 
        geothermal plants & heat pumps \\ 
        H2-electrolysis (fuel cell) & heat cogeneration from power plants \\ 
        lignite plants & ~ \\ \cline{2-2}
        lithium-ion batteries & \textbf{transportation} \\ \cline{2-2}
        nuclear plants & diesel (incl. synthetic diesel) \\ 
        oil plants & gasoline \\ 
        power-to-methane & natural gas (incl. synthetic methane) \\ 
        power-to-diesel & battery electric  \\ 
        pumped-hydro storage & plug-in hybrids \\ 
        run-of-river power plants & fuel cells \\
        photovoltaics & ~ \\ 
        transmission grids (220 kV + 380 kV) & ~ \\ 
        waste incineration plants & ~ \\ 
        wind (onshore and offshore) & ~ \\\hline 
    \end{tabular}
		\label{secmod:techs}
\end{table}

If not stated otherwise, we use four typical periods, consisting of 12 time steps with the length of one hour, thereby resulting in a total of 48 typical time steps. The limited temporal resolutions allows us to perform a benchmark computation using classical methods. We derived the typical periods using the TSAM package \citep{Hoffmann.2020}. For the time series aggregation, we use hierarchical clustering. As a benchmark, we solve the optimization problem at full spatial resolution without employing SpArta. 

We implemented the SpArta method and modeled the case study in the SecMOD framework \citep{Reinert.2022}, using Python $3.7.5$ and Pyomo $5.7.1$.  Computations were performed at the high-performance cluster of RWTH Aachen University (CLAIX-2018-MPI, 2 Intel Xeon Platinum 8160 Processors with 2.1 GHz, 4 GB main memory per core). As a solver, we use the Gurobi persistent algorithm with Gurobi $9.0.0$, and parallelization for up to $8$ cores. If not stated otherwise, we used the fast-forward algorithm to determine the spatial resolution during the iteration. 

\begin{figure}[ht]
	\centering
	\def\svgwidth{175pt}
	\includegraphics[width=\textwidth]{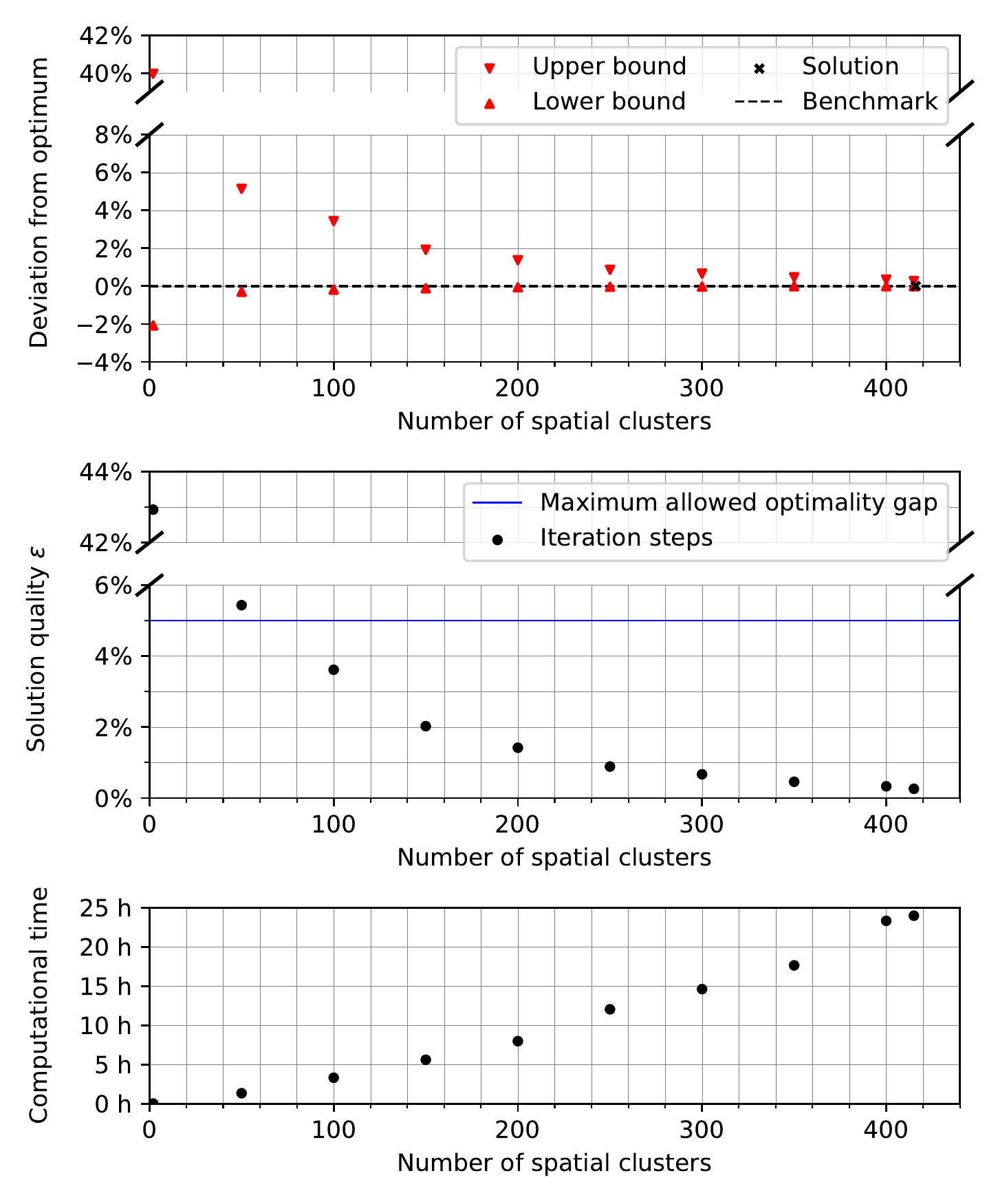}
	\caption[Convergence, solution quality and computational time.]{Convergence of upper and lower bound towards the solution of the full-scale problem with a continuously increasing spatial resolution at the step size of 50 nodes per iteration (top). The solution quality $\epsilon$ with respect to the lower bound is shown in the middle. Computational time is shown at the bottom.}	
	\label{fig:convergence}
\end{figure}

Figure \ref{fig:convergence} shows the convergence of the lower and upper bound with increasing spatial resolution when using SpArta. To show the convergence of lower and upper bounds, we did not define an optimality gap as a termination criterion but ran SpArta up to the full resolution of the benchmark with a continuous increase in spatial resolution. Here, the upper bound value shows the result of the upper bound relaxation of the aggregated system without local redesign relative to the benchmark. The lower and upper bound solutions converge towards the solution of the benchmark as the number of clusters increases.
\FloatBarrier

We observed no significant differences between the aggregation methods k-means, k-medoids, and the hierarchical algorithm in a preliminary investigation. Hence, the following results are obtained using the k-medoids algorithm. In our case study with $48$ typical time steps, a spatial resolution of slightly above $50$ nodes is necessary to obtain an optimal solution within the optimality gap of $5 \%$, i.e., around $15$ $\%$ of the original spatial resolution. Compared to the benchmark, the solution with SpArta only requires 17 \% of the computational time (given in wall-clock hours).

\begin{figure}[ht]
	\centering
	\includegraphics[width=\textwidth]{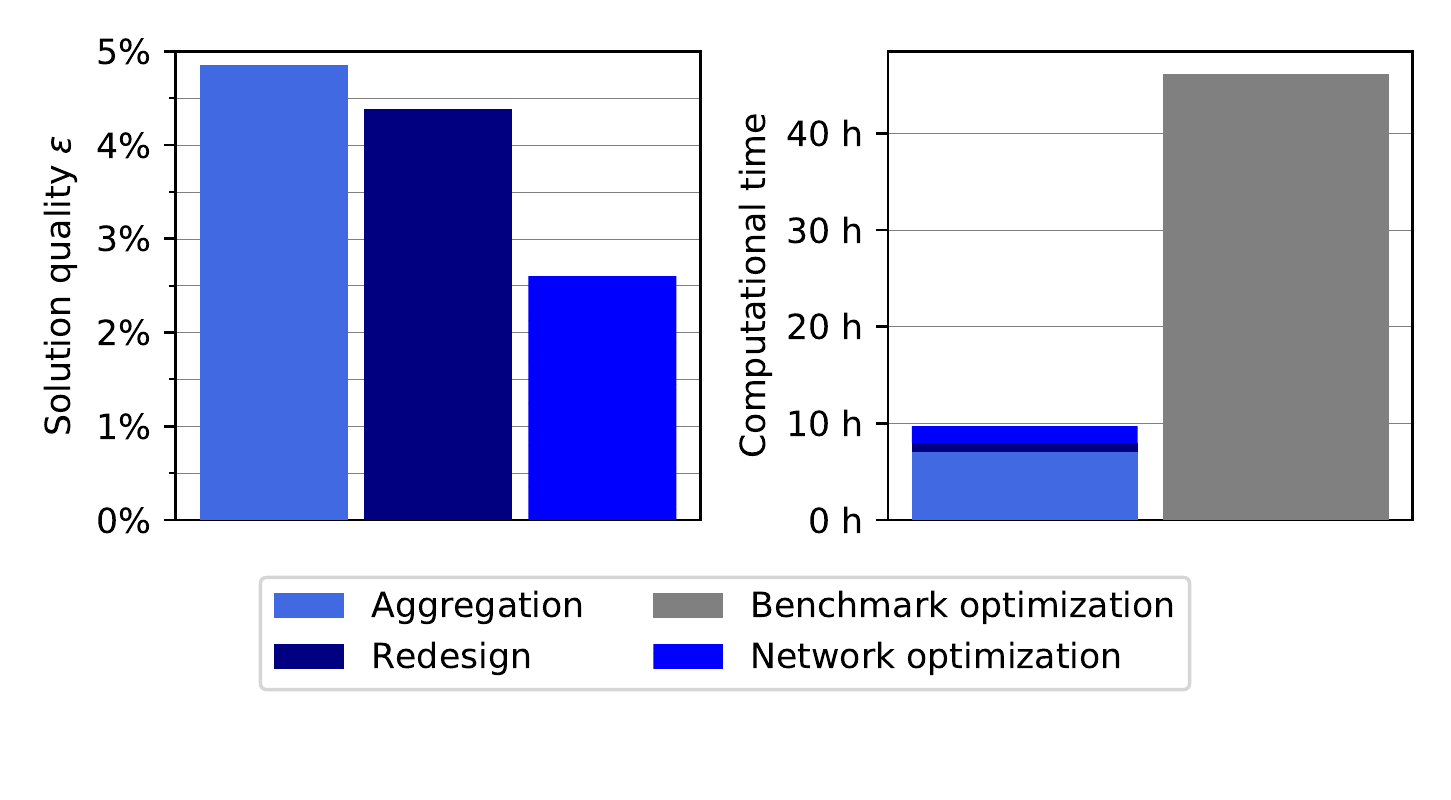}
	\caption[Solution quality and computational time for each step.]{Solution quality and computational time for finding the aggregated solution (aggregation), the decomposed solution at full scale (redesign), and the network optimization. The benchmark optimization is the optimal solution calculated with classical methods at full scale without SpArta.}	
	\label{fig:addopt}
\end{figure}

The actual solution quality of SpArta is higher than the predefined optimality gap $\epsilon_{SpArta}$ due to the subsequent local redesign and network optimization, which improve upon the upper bound. Figure \ref{fig:addopt} shows that with $\epsilon_{SpArta}= 5 \%$, the bounded error after the redesign is at less than $4.4$ $\%$. Further, with the network optimization, the bounded error is reduced to around $2.5$ $\%$. Hence, we observe a significant increase in solution quality after the decomposition to full spatial resolution. The benchmark optimization is the globally optimal solution calculated with classical methods at full-scale without SpArta and requires significantly higher computational time.

The local redesign to full spatial resolution consumes relatively little time compared to the iterative aggregated optimization (< 10 \% of the overall computational time). As the redesign problems are independent of each other, they can be parallelized, which has not even been fully exploited here. Subsequent to the local redesign, we performed an additional network optimization of the full-scale problem, i.e., an operational optimization including a design optimization of the grid, as discussed in Section \ref{subsub:networkopt}. The network optimization improves the solution quality significantly and further leads to a feasible solution with the DC-load-flow approach. Further, the computational effort is relatively low (< 18 \% of the computational time) compared to the gain in solution quality.

\FloatBarrier  

\begin{figure}[ht]
	\centering
	\includegraphics[width=\textwidth]{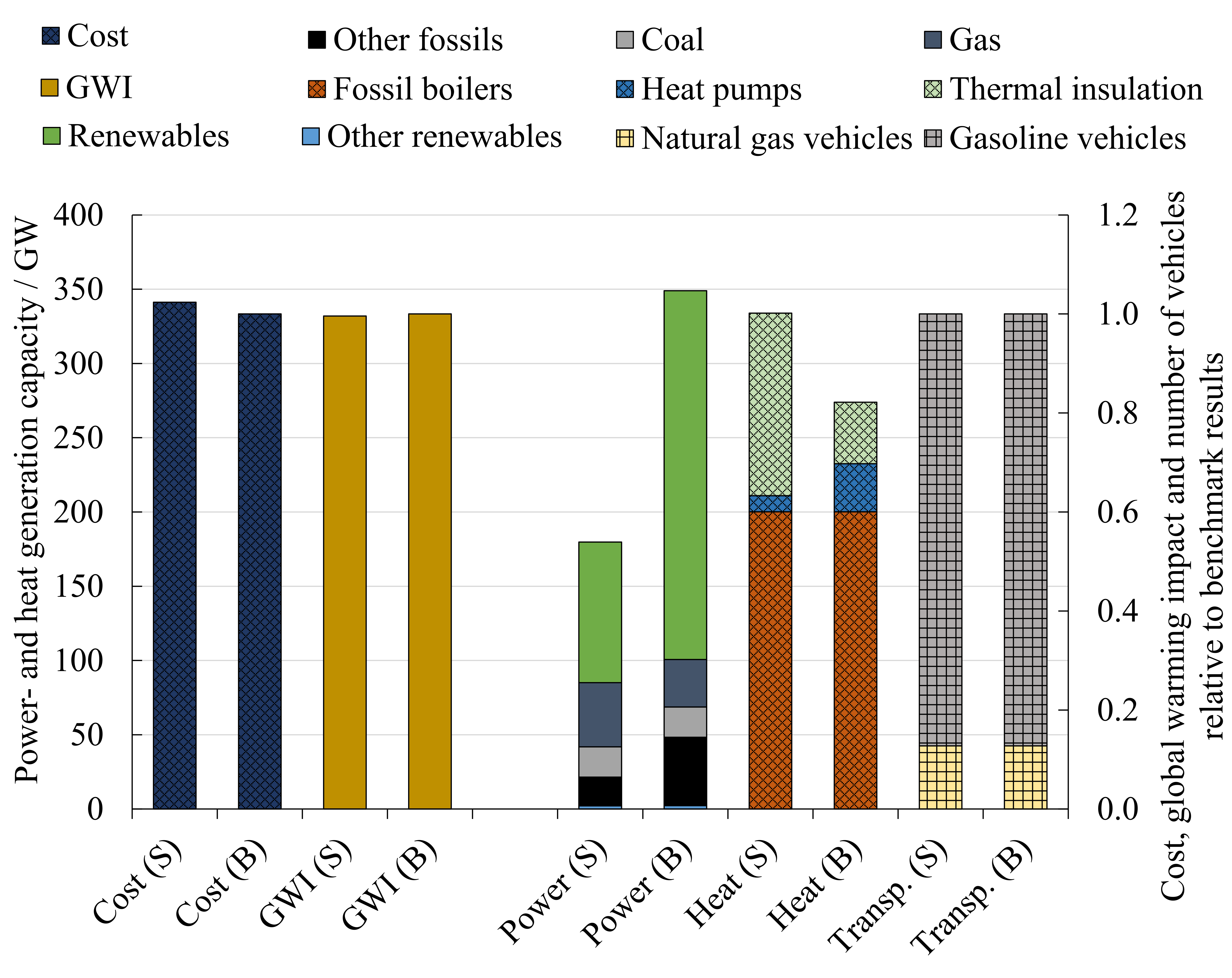}
	\caption{Comparison between the resulting system designs derived using SpArta (S) and the benchmark (B) regarding system cost, global warming impact (GWI) and capacity in the electricity-, heating-, and transport sectors.}	
	\label{fig:companalysis}
\end{figure}

Figure \ref{fig:companalysis} compares the resulting redesigned system using SpArta and the benchmark regarding cost, global warming impact, and further the electricity, heat, and transport capacity. We observe that while SpArta can identify a near-optimal solution that satisfies the environmental targets (here shown as global warming impact), the resulting system design differs from the benchmark design, computed at full spatial resolution. While the benchmark expands renewable converter capacity (wind and photovoltaics), SpArta leads to a system design with lower renewable energy converter capacity but larger shares of thermal insulation in the heating sector. 
This difference in capacity is due to the fact that SpArta uses the upper bound for the local redesign. Hence, weather conditions for fluctuating renewables are estimated pessimistically. The lower availability of renewable electricity shifts the preference from electricity-based heating via heat pumps to reduce the heat demand via thermal insulation. For this reason, SpArta design are near-optimal solutions that are conservative regarding the availability of renewables. The valid, near-optimal system design identified by SpArta is demonstrating the rich near-optimal design space in energy system modeling which should be systematically explored for decision support (as e.g., in \cite{Hennen.2017} and \cite{Lombardi.2020}). 

\FloatBarrier

\begin{figure}[ht]
	\centering
	\includegraphics[width=\textwidth]{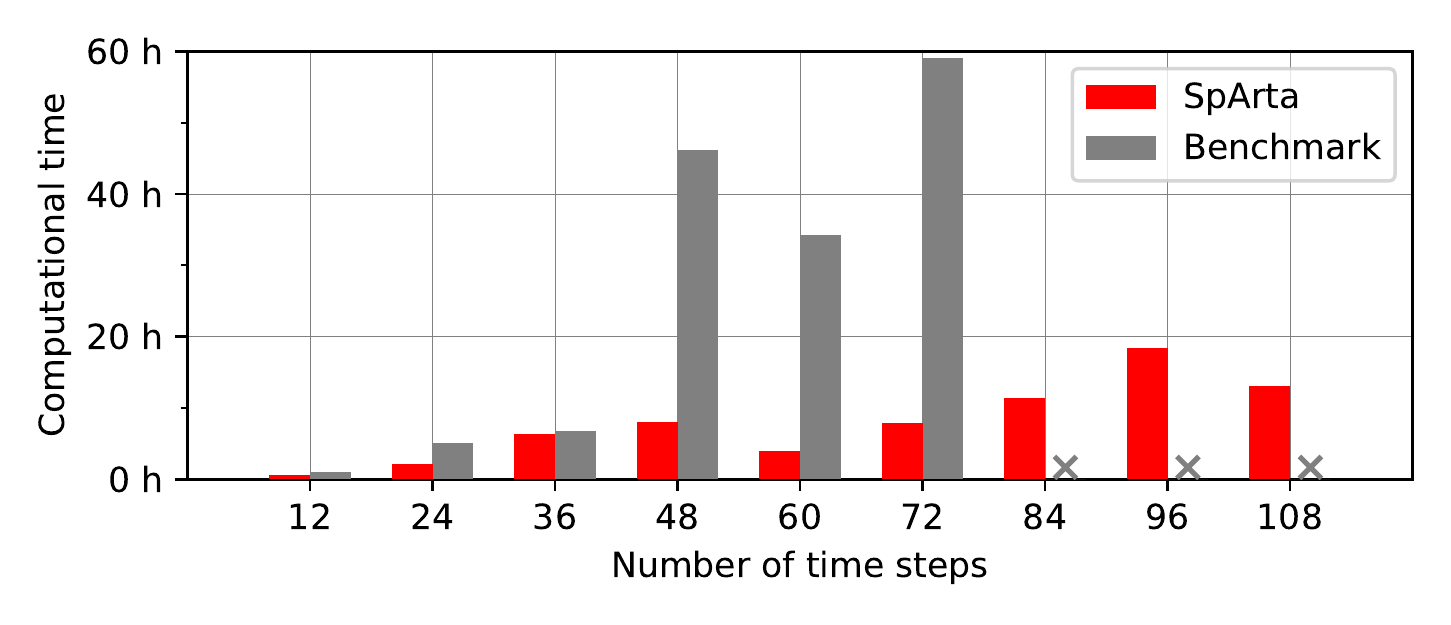}
	\caption[Computational time over problem size.]{Computational time over problem size, represented by the number of time steps. The required optimality gap $\epsilon_{SpArta}$ is $5 \%$. However, the final solution quality was better than 3 \% in all cases. If the problem could not be solved within the computational time limit of $72$ hours, we marked the problem with an 'x'.}	
	\label{fig:comptime}
\end{figure}

In Figure \ref{fig:comptime}, we compare the computational time of problems with different temporal aggregations with a time limit of $72$ hours. As a benchmark, we solve the original problem at full spatial resolution. Using SpArta, we optimize the system with a subsequent local redesign and a network optimization. We set the optimality gap $\epsilon_{SpArta}$ to a maximum of $5$ $\%$. However, the average final solution quality is around $2$ $\%$ and does not exceed $3$ $\%$ in any of the cases. 

For the largest problem size that could be solved for the benchmark on the available hardware, we observe a reduction in computational time by a factor of $7.5$ when SpArta is used. Further, we observe an increase in solvable problem size within the time limit by a factor of $1.5$ (Figure \ref{fig:comptime}). Hence, SpArta allows solving problems that could not be solved without SpArta within the limits of computational time. Furthermore, the computational time follows a quasi-linear curve: Employing the SpArta method at different problem sizes, we observe a significantly smaller increase in computational effort compared to the full-scale problem. In comparison to other state-of-the-art methods (e.g., \cite{Syranidou.2020, Reinert.2020}), SpArta finds a solution to the full-scale problem with known solution quality. 

\FloatBarrier	
\section{Conclusions}\label{sec:conclusion}
High spatial resolution is necessary to optimize energy systems with high shares of renewable energy. However, high spatial resolution simultaneously increases computational complexity, which leads to a trade-off between accuracy and computational time. In this work, we propose the SpArta method for rigorous optimization of regionally resolved energy systems by spatial aggregation and decomposition. SpArta significantly decreases computational effort while maintaining full spatial resolution. In SpArta, we find lower and upper bounds for the problem at lower spatial resolution. We increase spatial resolution iteratively, until the required optimality gap is reached. To obtain a solution at full spatial resolution, we solve the design problem by decomposing the full problem. We employ the results of the aggregated optimization for an additional design optimization of each cluster. Recombining the fully resolved clusters leads to a near-optimal design and operation at full spatial resolution.

Sparta allows to rigorously solve problems with larger problem size. Applying SpArta to regionally resolved energy system optimizations can significantly reduce computational effort while finding a feasible solution at high spatial resolution within the optimality gap: For a sector-coupled national energy system with $416$ nodes, we demonstrate an increase of solvable problem size by more than $50$ $\%$. Further, we observe that SpArta decreases computational time by a factor 7.5 in our case study and shows a quasilinear increase of computational time with problem size. The trend in computational time renders SpArta a promising method for energy systems with large problem size. Further parallelization could even amplify the benefits of SpArta. 

As the spatial extent of interconnected market zones and interacting energy systems increases, energy system models grow in problem size and complexity. The SpArta method enables the feasible synthesis of future sector-coupled energy systems at high spatial resolution with known solution quality.

\section*{Author contributions}
\textbf{Christiane Reinert}: Conceptualization, Methodology, Software, Investigation, Data Curation, Writing - Original Draft, Visualization, Project administration. 
\textbf{Benedikt Nilges}: Writing - Review \& Editing, Methodology, Software, Investigation, Visualization.
\textbf{Nils Baumgärtner}: Conceptualization, Methodology, Supervision, Writing - Review \& Editing. 
\textbf{André Bardow}: Conceptualization, Methodology, Resources, Writing - Review \& Editing, Supervision, Funding acquisition.

\section*{Declaration of Competing Interest}
We have no conflict of interest.

\section*{Acknowledgements}
Christiane Reinert acknowledges financial support by the Ministry of Economics, Innovation, Digitalization and Energy of North-Rhine Westphalia (Grant number: EFO 0001G). This work was partly sponsored by the Swiss Federal Office of Energy’s “SWEET” programme and performed in the “PATHFNDR“ consortium. Simulations were performed with computing resources granted by RWTH Aachen University under project rwth0755. The support is gratefully acknowledged. The authors thank Benedict Brosius and Julia Frohmann for their contributions to the case study.

\section*{Nomenclature}\label{app:Nomenclature}
Table \ref{tab:Nomenclature} gives a brief overview of all sets, superscripts, parameters, and variables used to describe the energy system and the SpArta approach. Furthermore it defines additional constraints to restrict parameters and variables to extend Section \ref{sub:GenericproblemDefinition}.
\begin{table}[ht]
	\caption{Nomenclature.}
	\label{tab:Nomenclature}
	\begin{tabularx}{\textwidth}{llll}

	\hline
	\multicolumn{2}{l}{\textbf{sets}} & \multicolumn{2}{l}{\textbf{superscripts}}\\
	$t\in T$ & time step & $k^{env}$ & environmental\\
	$n\in N$ & node & $k^{inv}$ & invest\\
	$c\in C$ & component & $k^{op}$ & operation \\
	$y\in Y$& investment year & $k^{eco}$ & economic\\
	$b\in B$ & product & $p^{flow/loss}$ & transport/loss\\
	$l\in L$ & edge & $p^{nom}$ & nominal\\
	$a\in A$ & spatial cluster & $p^{imp/exp}$ & import/export\\
	& &  $p^{min/max}$ & minimal/maximal\\
	& & $p^{int/ext}$ & internal/external\\
	& & $C^{grid}$ & grid\\
	\multicolumn{2}{l}{\textbf{parameters}} & $C^{sto}$& storage \\
	$k$ & specific cost & $C^{prod}$ & production\\
	$\beta$ & net present value factor &  & \\
	$d$ & exogenous demand & \multicolumn{2}{l}{\textbf{variables}} \\
	$h$ & time horizon&  $p$ & production\\
	$\eta$ & specific losses &&\\
	$\alpha$ & availability &&\\
	$\theta$ & product ratio matrix &&\\
	$\gamma$ & capacity factor &&\\
	$\sigma$ & connection between nodes and edges &&\\ 
	$\delta$ & needed secured capacity &&\\
	$\lambda$ & overall demand needs &&\\
	$\omega$ & share of useable capacity &&\\
	$\Delta t_t$& length of a time step \\
	\multicolumn{4}{l}{}\\
	\multicolumn{4}{l}{\textbf{additional constraints}}\\
	\multicolumn{2}{l}{$p_{c,n,y}^{nom} \in \mathbb{R}^+ \ ,\forall c\in C^{prod}, \forall n\in N, \forall y\in Y$} & \multicolumn{2}{l}{nominal capacity} \\
	\multicolumn{2}{l}{$p_{c,l,y}^{flow,nom} \in \mathbb{R}^+ \ ,\forall c\in C^{grid}, \forall l\in L, \forall y\in Y$} & \multicolumn{2}{l}{nominal grid capacity} \\
	\multicolumn{2}{l}{$p_{c,l,t}^{flow} \in \mathbb{R} \ ,\forall c\in C^{grid}, \forall l\in L, \forall t\in T, \forall y\in Y$} & \multicolumn{2}{l}{transport flow} \\
	\multicolumn{2}{l}{$p_{b,t}^{flow,loss} \in \mathbb{R}^+ \ ,\forall b\in B, \forall l\in L, \forall t\in T$} & \multicolumn{2}{l}{transport loss} \\
	\multicolumn{2}{l}{$\eta_{c,l} \in [0,1] \ ,\forall c\in C^{grid}, \forall l\in L$} & \multicolumn{2}{l}{transport efficiency} \\
	\multicolumn{2}{l}{$p_{c,n,t} \in \mathbb{R}^+ \ ,\forall c\in C^{prod}, \forall n\in N, \forall t\in T$} & \multicolumn{2}{l}{production} \\
	\multicolumn{2}{l}{$\alpha_{c,n,t} \in [0,1] \ ,\forall c\in C^{prod}, \forall n\in N, \forall t\in T$} & \multicolumn{2}{l}{production availability} \\
	\multicolumn{2}{l}{$p_{b,n,t}^{exp} \in \mathbb{R} \ ,\forall b\in B, \forall c\in C^{prod}, \forall n\in N, \forall t\in T$} & \multicolumn{2}{l}{nodal product export} \\
	\multicolumn{2}{l}{$p_{b,t}^{import} \in \mathbb{R}^+ \ ,\forall b \in B, \forall t\in T$} & \multicolumn{2}{l}{global product import} \\
	\multicolumn{2}{l}{$d_{b,n,t} \in \mathbb{R}^+ \ ,\forall b\in B, \forall n\in N, \forall t\in T$} & \multicolumn{2}{l}{demand} \\
	\multicolumn{2}{l}{$\sigma_{l,n,m} \in \{0,1\} \ ,\forall l\in L, \forall m,n \in N, \forall m\in N\backslash\{n\}$} & \multicolumn{2}{l}{connectivity matrix} \\
	\multicolumn{2}{l}{$\omega_{c,n,t} \in [0,1] \ ,\forall c\in C^{prod}, \forall n \in N,\forall t \in T$} & \multicolumn{2}{l}{share of useable capacity}\\
	\multicolumn{2}{l}{$\theta_{b,c} \in \mathbb{R} \ ,\forall b\in B, \forall c\in C$} & \multicolumn{2}{l}{product ratio matrix} \\
	\multicolumn{2}{l}{$e^{max}, \beta, \gamma, k, p^{min}, h_c \in \mathbb{R}^+$} & \multicolumn{2}{l}{emission limit, }\\
	& & \multicolumn{2}{l}{net present value factor, } \\
	& & \multicolumn{2}{l}{capacity factor, costs, }\\
	& & \multicolumn{2}{l}{minimal capacity limit,} \\
	& & \multicolumn{2}{l}{and lifetime} \\
	\hline
\label{tab:Nomenclature}
\end{tabularx}
\end{table}
\FloatBarrier

  \bibliographystyle{apalike}
  \renewcommand{\refname}{Bibliography}  
  \bibliography{literature.bib}

\end{document}